\documentclass[11pt]{article}

\usepackage[margin=1in]{geometry}
\usepackage{xcolor}
\usepackage{amsmath,amssymb,amsthm,mathtools}
\usepackage{stmaryrd}
\usepackage{enumitem}
\usepackage{float}
\usepackage{tabularx}
\usepackage{multirow}
\usepackage{booktabs}
\usepackage{hyperref}
\usepackage[nameinlink,capitalise,noabbrev]{cleveref}

\hypersetup{
  colorlinks=true,
  linkcolor=blue!50!black,
  citecolor=blue!50!black,
  urlcolor=blue!50!black
}

\newtheorem{theorem}{Theorem}[section]
\newtheorem{proposition}[theorem]{Proposition}
\newtheorem{lemma}[theorem]{Lemma}

\theoremstyle{definition}
\newtheorem{definition}[theorem]{Definition}

\newtheorem{question}[theorem]{Question}

\newcommand{\rank}{\operatorname{rank}}
\newcommand{\spt}{\operatorname{spt}}

\newcommand{\R}{\mathbb R}
\newcommand{\Z}{\mathbb Z}

\newcommand{\cR}{\mathcal R}

\DeclareMathOperator{\Gr}{Gr}

\DeclareMathOperator{\Sing}{Sing}
\newcommand{\MOD}{\mathbf{MOD}}

\newcommand{\cH}{\mathcal H}

\newcommand{\cC}{\mathcal C}
\newcommand{\bM}{\mathbf M}
\newcommand{\ip}[2]{\left\langle #1,#2\right\rangle}
\newcommand{\abs}[1]{\left|#1\right|}
\newcommand{\norm}[1]{\left\lVert#1\right\rVert}

\title{\textbf{Real Determinantal and Segre Varieties: Stability and Area Minimization}}

\author{Hongbin Cui\and Xiaoxiang Jiao \and Hangyue Zhu}
\date{\today}

\begin{document}

\maketitle

\begin{abstract}
We apply Liu's recently introduced moderation criterion to complete the classification of area-minimizing real determinantal varieties and real Segre varieties, thereby settling the remaining cases left open by Pedersen and Kerckhove--Lawlor \cite{Pedersen93,KerckhoveLawlor1999}. In particular, we prove that $C(2,3,1)$ is area-minimizing modulo $2$, resolving the $g=4$, $(m_1,m_2)=(1,1)$ case of the isoparametric focal-cone problem raised by Tang and Zhang \cite{TZ20} and recently revisited by Liu \cite{Liu2026}.
\end{abstract}

\noindent\textbf{Keywords.} Area-minimizing; Determinantal varieties; Segre varieties; moderation.

\noindent\textbf{2020 Mathematics Subject Classification.} 49Q05, 49Q15, 53A10, 14M12.
\tableofcontents
\section{Introduction}

Determining the stability and area-minimizing properties of minimal cones is a central problem in geometric measure theory. Over the past several decades, several powerful tools have been developed to address these questions, including calibrations, Lawlor's curvature criterion, and Liu's recently developed moderation criterion \cite{HL82,angle,Lawlorcriterion,Liu2026}. These tools have been successfully applied to many fundamental examples, including the Simons cone, unions of planes, and cones over various homogeneous spaces and their products \cite{Simons1,Nance,angle,Kan,XU2018739,TZ20,JCX22}\footnote{Although the work of Bombieri--De Giorgi--Giusti \cite{Simons1} predates the calibration method in \cite{HL82}, their methods are later reinterpreted in terms of calibrations; see \cite{Simons1.5,Simonscone2}.}.

Two of the most natural families in which to study these problems are the real determinantal cones  and the Segre cones:

\begin{definition}[Real determinantal cones (varieties)]
  For $1\leq r<p\leq q$, the \emph{real determinantal cone} is defined by\[
C(p,q,r):=\{A\in M_{p\times q}(\mathbb R):\operatorname{rank}A\leq r\}.
\]The cone has dimension $r(p+q-r)$ \cite[Corollary 4.3]{KerckhoveLawlor1999}. Its regular points are exactly the matrices of rank $r$, while $\operatorname{Sing}(C(p,q,r))=C(p,q,r-1)$ \cite[Corollary 4.4]{KerckhoveLawlor1999}. Setting $C_r=C(p,q,r)$, we have:\[
\{0\}=C_0\subset C_1\subset\cdots\subset C_{p-1},\qquad\operatorname{Sing}(C_r)=C_{r-1}.
  \] 
We also define its deficit by $\delta:=p+q-2r$. Then $\operatorname{codim}_{C(p,q,r)}\operatorname{Sing}(C(p,q,r))=\delta+1$.
\end{definition}

\begin{definition}[Segre cones (varieties)]For integers $d\ge 2$ and $n_1,\cdots,n_d\ge 2$, the \emph{Segre cones} are defined by \[
  C(\mathbb S^{n_1-1}\otimes\cdots\otimes\mathbb S^{n_d-1})\hookrightarrow\mathbb R^{\,n_1\cdots n_d}.
  \]
The cone has dimension $1-d+\sum_{i=1}^d n_i$, and its unique singular point is the origin. 
\end{definition}

The two families are closely related by the following identification:\begin{equation}\label{eq:identicication}
  C(\mathbb S^{p-1}\otimes \mathbb S^{q-1})=\{A\in M_{p\times q}(\mathbb R):\operatorname{rank}A\le1\}=C(p,q,1).
\end{equation}
Thus a two-factor Segre cone is precisely a rank-one determinantal cone. Earlier works determine their minimizing behavior for most parameter values, leaving only a small number of low-deficit determinantal cones and low-dimensional Segre cones outside the scope \cite{Pedersen93,KerckhoveLawlor1999,CUI26}. These remaining cases naturally lead to the following two questions:

\begin{question}\label{quest:det}
What is the minimizing behavior of the exceptional determinantal cones of \cite{KerckhoveLawlor1999}?
\end{question}
\begin{question}\label{quest:segre}
What is the minimizing behavior of the exceptional Segre cones of \cite{Pedersen93,CUI26}?
\end{question}

The status of Questions \ref{quest:det} and \ref{quest:segre} before and after the present work are summarized as in the following Table \ref{tab:introduction-status}. Throughout the paper, minimization is understood with integer coefficients when the link is orientable and modulo 2 otherwise.

\begin{table}[H]
\centering
\small
\renewcommand{\arraystretch}{1.25}
\begin{tabularx}{\textwidth}{
  >{\raggedright\arraybackslash}p{0.27\textwidth}
  >{\raggedright\arraybackslash}p{0.31\textwidth}
  >{\raggedright\arraybackslash}X}
\toprule
Family & Previously known status & Result established in this paper \\
\midrule
$C(p,q,r)$, $\delta\geq4$
& Classified by \cite{Pedersen93,KerckhoveLawlor1999}
& N/A \\

$C(p,p,p-1)$, $\delta=2$
& Non-minimizing by \cite{KerckhoveLawlor1999}
& Unstable for every $p\geq2$ \\

$C(r+1,r+2,r)$, $\delta=3$
& Exceptional case of \cite{Pedersen93,KerckhoveLawlor1999}
& Mod $2$ area-minimizing \\

Two-factor Segre cones
& Rank-one determinantal cones
& Classified by determinantal results \\

Segre cones with $d\geq3$, dimension $n\ge 5$
& Area-minimizing by \cite{CUI26}
& N/A \\

$C(\mathbb S^1\otimes\mathbb S^1\otimes\mathbb S^1)$
& Exceptional case of \cite{CUI26}
& Area-minimizing \\
\bottomrule
\end{tabularx}
\caption{Status of the exceptional determinantal and Segre cones.}
\label{tab:introduction-status}
\end{table}

\paragraph{Main results.}
Our first main result is that every determinantal cone $C(r+1,r+2,r)$, $r\geq1$, is mod $2$ area-minimizing. We prove this by applying Liu's moderation criterion. In particular, the case $C(2,3,1)$, which is isomorphic to $C(\mathbb S^1\otimes\mathbb S^2)$, is mod $2$ area-minimizing and answers the question in \cite[Question 15]{Pedersen93}.

The same cone  $C(2,3,1)$ also has a natural interpretation in isoparametric geometry. Tang and Zhang \cite{TZ20} singled out the low-multiplicity cases $(m_1,m_2)=(1,1)$, for $g=3,4,6$, as the cases not covered by the known area-minimization results for cones over focal submanifolds, and asked whether the two focal cones $C(M_{+})$ and $C(M_{-})$ in each of these cases are area-minimizing. The two cones for the $g=3$ case are both the cone over the Veronese minimal embedding of $\mathbb{R}P^2$, which was recently proved by Liu \cite{Liu2026} to be area-minimizing modulo $2$. Liu also further suggested that the moderation method might also apply to isoparametric focal cones.

For $g=4$ and $(m_1,m_2)=(1,1)$, the two focal submanifolds admit explicit geometric descriptions; see, for example \cite[pp.~155--159]{CR15}. With the convention of \cite{TZ20}, $M_+$ is the Stiefel manifold $V_2(\mathbb{R}^3)$ of orthonormal two-frames in $\mathbb{R}^3$. The cone $C(M_+)$ is in fact calibrated and thus area-minimizing;  let's give a brief explanation here (also see \cite[Sect. 1.2]{CuiFKM}). Under the identification $\mathbb{R}^3\oplus\mathbb{R}^3\cong\mathbb{C}^3$, $C(M_+)$ can be written as the complex quadric
\[C(M_+)=\left\{z\in\mathbb{C}^3:z_1^2+z_2^2+z_3^2=0\right\}.
\]
It is therefore calibrated by the power of Kähler forms. The other focal submanifold is
\[M_-=\left\{e^{\mathrm{i}\theta}x:\theta\in\mathbb{R},\ x\in S^2\right\}\cong(S^1\times S^2)/\{(u,x)\sim(-u,-x)\},
\]
and its cone is precisely
\[C(M_-)=C(S^1\otimes S^2)=C(2,3,1).\]
Thus, the modulo $2$ area-minimization of $C(2,3,1)$ established in the present paper completes the $g=4$, $(1,1)$ case. Consequently, the only multiplicity case in the question of Tang and Zhang that remains open is $g=6$, $m=1$: are both focal cones $C(M_+)$ and $C(M_-)$ area-minimizing modulo $2$? Here $M_\pm^5\subset S^7$ are noncongruent and both diffeomorphic to $S^3\times\mathbb{R}P^2$; see \cite{Miy93} for their geometry.

Our second result is that $C(p,p,p-1)$ is unstable for all $p\ge 2$. This strengthens the previously known non-minimizing result of Kerckhove--Lawlor \cite[Theorem 8.3(1)]{KerckhoveLawlor1999}.

Finally, we construct an explicit real orthogonal isomorphism between the triple Segre cone $C(\mathbb S^1\otimes\mathbb S^1\otimes\mathbb S^1)$ and the Harvey--Lawson special Lagrangian cone in $\mathbb C^4$. Consequently, $C(\mathbb S^1\otimes\mathbb S^1\otimes\mathbb S^1)$ is calibrated and area-minimizing.

\paragraph{Organization.} Section~\ref{sec:moderation} recalls Liu's moderation criterion. Section~\ref{sec:remaining-real-determinantal} proves the mod $2$ area minimization of the family $C(r+1,r+2,r)$. Section~\ref{sec:unstability-Cppp-1} proves the instability of the family $C(p,p,p-1)$. Combining these two sections completes the full classification of determinantal cones. Finally, Section~\ref{sec:Segre} identifies the triple tensor Segre cone with the Harvey--Lawson cone and obtain the complete classification of Segre cones.

\section{Liu's moderation criterion}\label{sec:moderation}

Throughout the paper, we use the standard notation and basic results for mod $2$ currents, slices, push-forwards, and the area formula, as in \cite{FedGMT}. Moreover, if $P,Q$ are unoriented $d$-planes, then $\abs{\ip{P}{Q}}$ denotes the absolute value of the inner product of unit simple $d$-vectors representing $P,Q$.

We first record the area-minimizing criterion introduced by Liu \cite[Definition~3.0.1 and Theorem~4]{Liu2026}, which was inspired by Morgan's earlier results \cite{Morgan1,MORGAN2}.

\begin{definition}[Moderation; see {\cite[Definition 3.0.1]{Liu2026}}]\label{def:moderation} Let $T\in\mathbf I_d(\R^{d+c};\Z_2)$ be a compactly supported current. For $W\in\Gr_d(\R^{d+c})$, define
\begin{equation}
 I_T(W):=\int \abs{\ip{W}{\operatorname{Tan}(T,x)}}\,d\norm{T}(x),
\label{eq:projection-functional}
\end{equation}
where $\operatorname{Tan}(T,x)$ is the approximate tangent plane of $T$ at $x$. We say that $T$ is a \emph{moderation} if, for $\norm{T}$-almost every $x$, with $E_x=\operatorname{Tan}(T,x)$,
\begin{align}
I_T(E_x)&=\sup_{W\in\Gr_d(\R^{d+c})}I_T(W),\label{eq:moderation-one}\\
\bM\bigl((\pi_{E_x})_\#T\bigr)&=I_T(E_x)>0.\label{eq:moderation-two}
\end{align}
The common value in \eqref{eq:moderation-one} and \eqref{eq:moderation-two} is denoted by
\[
\MOD(T):=\sup_{W\in\Gr_d(\R^{d+c})}I_T(W).
\]
\end{definition}

Condition \eqref{eq:moderation-one} says that almost every tangent plane of $T$ maximizes the averaged projection functional.  Condition \eqref{eq:moderation-two} says that projection onto such a tangent plane produces no cancellation; its mass is exactly the integral of the projected Jacobian.

\begin{theorem}[Liu's moderation criterion; see {\cite[Theorem 4]{Liu2026}}]
\label{thm:liu}
Let $T\in\mathbf I_d(\R^{d+c};\Z_2)$ be a moderation. Let $S\in\mathbf I_d(\R^{d+c};\Z_2)$ and suppose that
\begin{equation}S-T=\partial Q\label{eq:homological-comparison}
\end{equation}
for some $Q\in\mathbf I_{d+1}(\R^{d+c};\Z_2)$. Then $\bM(T)\leq \bM(S).$ 
\end{theorem}
\begin{proof}
For every $W\in\Gr_d(\R^{d+c})$, \eqref{eq:homological-comparison} gives
\[
(\pi_W)_\#S-(\pi_W)_\#T=\partial\bigl((\pi_W)_\#Q\bigr).
\]
Since the image of $\pi_W$ has dimension $d$, its $(d+1)$-Jacobian vanishes. Hence the push-forward estimate gives
\[
\bM\bigl((\pi_W)_\#Q\bigr)\leq \int J_{d+1}\pi_W\,d\norm{Q}=0,
\]
so that $(\pi_W)_\#Q=0$ and therefore
\begin{equation}
(\pi_W)_\#S=(\pi_W)_\#T\qquad\text{for every }W\in\Gr_d(\R^{d+c}).
\label{eq:projected-currents-equal}
\end{equation}

For $\norm{T}$-almost every $x$, set $E_x=\operatorname{Tan}(T,x)$. By \eqref{eq:moderation-two} and \eqref{eq:projected-currents-equal}, we have the following inequality:
\[
\MOD(T)=\bM\bigl((\pi_{E_x})_\#S\bigr)\leq \int\abs{\ip{E_x}{\operatorname{Tan}(S,y)}}\,d\norm{S}(y).
\]
Integrating in $x$ and applying Tonelli's theorem, we obtain
\begin{align*}
\MOD(T)\bM(T)
&\leq \int\!\int\abs{\ip{\operatorname{Tan}(T,x)}{\operatorname{Tan}(S,y)}}\,d\norm{S}(y)\,d\norm{T}(x) \\
&=\int I_T\bigl(\operatorname{Tan}(S,y)\bigr)\,d\norm{S}(y) \\
&\leq \MOD(T)\bM(S),
\end{align*}
where the last inequality follows from \eqref{eq:moderation-one}. Dividing by $\MOD(T)>0$ proves the theorem.
\end{proof}

\section{The real determinantal cone \texorpdfstring{$C(r+1,r+2,r)$}{C(r+1,r+2,r)}}\label{sec:remaining-real-determinantal}

In this section, we use Liu's moderation criterion to prove that the determinantal cones $C(r+1,r+2,r)$ are mod 2-minimizing. This completes the classification of real determinantal cone as in the following Table \ref{tab:real-determinantal-cones}:

\begin{table}[H]
\centering
\small
\renewcommand{\arraystretch}{1.25}
\begin{tabularx}{\textwidth}{
  >{\centering\arraybackslash}p{0.12\textwidth}
  >{\centering\arraybackslash}p{0.25\textwidth}
  >{\centering\arraybackslash}p{0.16\textwidth}
  >{\raggedright\arraybackslash}X}
\toprule
$p+q-2r$ & Real determinantal cone & Dimension & Status and Source \\
\midrule
$2$
& $C(p,p,p-1)$
& $p^2-1$
& Non-minimizing by \cite[Theorem 8.3(1)]{KerckhoveLawlor1999}; in fact, unstable by Section~\ref{sec:unstability-Cppp-1}. \\
\addlinespace[0.6ex]
$3$
& $C(p,p+1,p-1)$
& $(p-1)(p+2)$
& Mod 2-minimizing by Section~\ref{sec:remaining-real-determinantal}. \\
\addlinespace[0.6ex]
\multirow{2}{*}{$\ge4$}
& $C(p,q,r)$, $p+q$ even
& $r(p+q-r)$
& Area-minimizing by \cite[Theorem 8.3(2)]{KerckhoveLawlor1999}. \\
\addlinespace[0.6ex]
& $C(p,q,r)$, $p+q$ odd
& $r(p+q-r)$
& Mod 2-minimizing by \cite[Theorem 8.3(3)]{KerckhoveLawlor1999}. \\
\bottomrule
\end{tabularx}
\caption{Classification of area-minimizing real determinantal cones. Here $1\le r<p\le q$.}
\label{tab:real-determinantal-cones}
\end{table}

We first determine the basic geometry of the cones $C(r+1,r+2,r)$. Set
\[p=r+1,\qquad q=r+2,\qquad V=M_{p,q}(\R),\qquad d=r(r+3)=pq-2,
\]
and define the cone and its regular stratum by
\[
\cC_r=C(r+1,r+2,r),\qquad\cR_r=\{X\in V:\rank X=r\}.
\]
One would have:
\begin{equation}
\dim\cC_r=r(p+q-r)=r(r+3)=d.\label{eq:dimension-main-stratum}
\end{equation}
The singular set of $\cC_r$ is the next lower-rank variety:
\[
\Sing(\cC_r)=C(r+1,r+2,r-1),
\]
and therefore
\begin{equation}
\dim\Sing(\cC_r)=(r-1)(p+q-r+1)=(r-1)(r+4)=d-4.
\label{eq:sing-dimension}
\end{equation}
In particular, the singular set is $\cH^{d-1}$-null.  Let
\[
C_r^\infty:=\llbracket\cC_r\rrbracket_2,\qquad T_r:=C_r^\infty\mathbin{\llcorner}\mathbb B_1(0).
\]
Thus $C_r^\infty$ is the multiplicity-one locally integral current modulo $2$ carried by the full cone, while $T_r$ is its unit-ball truncation. Since $\cR_r$ is a smooth manifold without boundary which carries the full $\cH^{d}$-measure, together with \eqref{eq:sing-dimension}, we have $\partial C_r^\infty=0$.

Let $\rho(X)=\abs{X}$. The restriction formula for currents gives
\[
\partial T_r=\partial\bigl(C_r^\infty\mathbin{\llcorner}\mathbb B_1(0)\bigr)=\langle C_r^\infty,\rho,1\rangle,
\]
where $\langle\,\cdot,\cdot,\cdot\,\rangle$ denotes slicing. Since $\cC_r$ is a cone, the level sets meet $\cR_r$ transversely, and therefore $\partial T_r=\llbracket\cC_r\cap\mathbb S^{pq-1}\rrbracket_2.$

Finally, since $\Sing(\cC_r)$ is $\cH^d$-null set, the approximate tangent plane of $T_r$ agrees with $T_X\cR_r$ for $\norm{T_r}$-almost every $X$.

\subsection{The normal Gauss map}

\begin{lemma}
\label{lem:tangent-normal}
For $X\in\cR_r$, set
\[
L_X=\ker X^T,\qquad K_X=\ker X.
\]
Then $\dim L_X=1$ and $\dim K_X=2$. If $v$ is a unit vector spanning $L_X$, then
\begin{align}
T_X\cR_r
&=\{Y\in V:v^TYk=0\text{ for every }k\in K_X\},\label{eq:tangent-formula}\\
N_X\cR_r
&=\{vk^T:k\in K_X\}=v\otimes K_X.
\label{eq:normal-formula}
\end{align}
\end{lemma}

\begin{proof}
Since $p=r+1$, $q=r+2$, and $\rank X=r$, the dimensional assertions of $L_X$ and $K_X$ are immediate. Choose orthogonal coordinates in the source and target so that
\[
X=\left(\begin{array}{c|c}\Sigma_{r\times r}&0_{r\times2}\\\hline 0_{1\times r}&0_{1\times2}\end{array}\right).
\]
where $\Sigma\in M_{r,r}(\R)$ is invertible. In a neighborhood of $X$, let $\gamma(t)\in\mathcal R_r$ be a smooth curve with $\gamma(0)=X$ and $\gamma'(0)=Y$. Then:\[\gamma(t)=
\begin{pmatrix}
\Sigma+tA+O(t^2)&tB+O(t^2)\\
tC+O(t^2)&tD+O(t^2)
\end{pmatrix}.\] The $\rank r $ condition is described by the Schur complement:
\[
\operatorname{rank}\gamma(t)=r\Longleftrightarrow tD+O(t^2)-\bigl(tC+O(t^2)\bigr)\bigl(\Sigma+tA+O(t^2)\bigr)^{-1}\bigl(tB+O(t^2)\bigr)=0.
  \]
Consequently, the tangent space consists precisely of those $Y$ whose lower-right $1\times2$ block vanishes. Equivalently,
\[
T_X\cR_r=\left\{\begin{pmatrix}
A&B\\
C&0\end{pmatrix}\right\}
=\{Y\in\operatorname{Hom}(\R^q,\R^p):Y(K_X)\subset v^\perp=\operatorname{Im}X\}.
\]
This proves \eqref{eq:tangent-formula}. Moreover, using $\ip{\cdot}{\cdot}_F$ to denote the Frobenius inner product $\langle A,B\rangle_F=\operatorname{tr}(A^TB)$, one has: $\langle Y,vk^T\rangle_F=\operatorname{tr}(Y^Tvk^T)=\operatorname{tr}(k^TY^Tv)=v^TYk=0.$ It follows that $v\otimes K_X$ is orthogonal to $T_X\cR_r$. Both spaces have dimension two, this proves \eqref{eq:normal-formula}.
\end{proof}
Since the moderation functional depends only on the distribution of tangent planes, it is convenient to work with their orthogonal complements. The normal plane at $X$ depends only on the line $[v]\in\mathbb{RP}^{p-1}$ and the two-plane $K_X\in\Gr_2(\R^q)$. We therefore define the Gauss map
\begin{equation}
\gamma:\cR_r\longrightarrow\mathbb{RP}^{p-1}\times\Gr_2(\R^q),\qquad\gamma(X)=([\ker X^T],\ker X).
\label{eq:gauss-map}
\end{equation}

\begin{lemma}
\label{lem:gauss-surjective}
The Gauss map \eqref{eq:gauss-map} is surjective. Its restriction to $\cR_r\cap \mathbb B_1(0)$ is also surjective.
\end{lemma}

\begin{proof}
Fix $([v],K)\in\mathbb{RP}^{p-1}\times\Gr_2(\R^q)$. Since $\dim K^\perp=\dim v^\perp=r$, choose a linear isomorphism $L:K^\perp\rightarrow v^\perp$ and define
\[
X=L\circ P_{K^\perp}:\R^q\longrightarrow\R^p,
\]
where $P_{K^\perp}$ is the orthogonal projection onto $K^\perp$. Then
\[
\rank X=r,\qquad \ker X=K,\qquad \ker X^T=(\operatorname{Im}X)^\perp=\R v.
\]
This proves surjectivity. Multiplying $X$ by a sufficiently small constant places it in $\mathbb B_1(0)$ without changing either kernel.
\end{proof}

Let $\sigma_p$ be the normalized $O(p)$-invariant probability measure on $\mathbb{RP}^{p-1}$, and let $\nu_q$ be the normalized $O(q)$-invariant probability measure on $\Gr_2(\R^q)$.

\begin{lemma}
\label{lem:gauss-haar}
Let $m_r=\bM(T_r)$, then
\begin{equation}
\gamma_\#\left(\cH^d\mathbin{\llcorner}(\cR_r\cap \mathbb B_1(0))\right)=m_r\,\sigma_p\otimes\nu_q.
\label{eq:gauss-haar}
\end{equation}
\end{lemma}

\begin{proof}
By \eqref{eq:sing-dimension}, $\cR_r$ has full $d$-dimensional measure in $\cC_r$. Since $\cC_r$ is a cone, its intersection with the unit sphere has dimension at most $d-1$. Consequently,
\[
\cH^d(\cR_r\cap \mathbb B_1(0))=\bM(T_r)=m_r<\infty.
\]

The group $O(p)\times O(q)$ acts orthogonally on $V$ by $(U,V_0)\cdot X=UXV_0^T.$ This action preserves $\cC_r$, $\cR_r$, $\mathbb B_1(0)$, and $\cH^d$. It is compatible with the Gauss map \eqref{eq:gauss-map} because
\[
\ker(UXV_0^T)^T=U(\ker X^T),\qquad\ker(UXV_0^T)=V_0(\ker X).
\]
Thus the pushforward measure in \eqref{eq:gauss-haar} is a finite $O(p)\times O(q)$-invariant measure on $\mathbb{RP}^{p-1}\times\Gr_2(\R^q).$ The product action is transitive, and $\sigma_p\otimes\nu_q$ is its unique invariant probability measure. Therefore the pushforward is a constant multiple of this product measure, using the definition of $m_r$ gives \eqref{eq:gauss-haar}. 
\end{proof}

\subsection{Verification of the moderation conditions}

The preceding Lemma \ref{lem:gauss-haar} converts integration over the determinantal cone into integration over the homogeneous space $\mathbb{RP}^{p-1}\times\Gr_2(\R^q).$ We now use this to verify moderation conditions in definition \ref{def:moderation}. We first pass between tangent and normal planes.

\begin{lemma}
\label{lem:hodge}
Let $E,W$ be unoriented $d$-planes in the $(d+2)$-dimensional Euclidean space $V$. Then
\[
\abs{\ip{E}{W}}=\abs{\ip{E^\perp}{W^\perp}}.
\]
\end{lemma}

\begin{proof}
Fix an orientation of $V$. For an unoriented plane $E$, we have $\xi_{E^\perp}=\pm *\xi_E$, where $*$ denotes the Hodge star. Since the Hodge star is an isometry, the two inner products agree up to sign.  Taking absolute values removes the dependence on orientations or signs.
\end{proof}

For $K\in\Gr_2(\R^q)$, the line $\Lambda^2K$ has two unit simple generators differing by a sign. Thus expressions of the form $\abs{\ip{\eta}{\xi_K}}$ do not depend on this choice, here $\xi_K$ denotes either generator. Define
\begin{equation}
c_q:=\int_{\Gr_2(\R^q)}\norm{\operatorname{proj}_{\Lambda^2K}(e_1\wedge e_2)}\,d\nu_q(K)=\int_{\Gr_2(\R^q)}\abs{\ip{e_1\wedge e_2}{\xi_K}}\,d\nu_q(K)>0.
\label{eq:def-cq}
\end{equation}
The integrand equals $1$ at $K=\operatorname{span}\{e_1,e_2\}$ and is positive on a nonempty open neighborhood of that plane, proving the strict positivity in \eqref{eq:def-cq}.

Equivalently, if $\theta_1,\theta_2$ denote the principal angles between $\operatorname{span}\{e_1,e_2\}$ and $K$, then $\bigl|\langle e_1\wedge e_2,\xi_K\rangle\bigr|=|\cos\theta_1\cos\theta_2|,$ and $c_q=\mathbb E_{\nu_q}[|\cos\theta_1\cos\theta_2|].$

\begin{proposition}
\label{prop:M1}
For every $W\in\Gr_d(V)$,
\begin{equation}
I_{T_r}(W)\leq \frac{m_rc_q}{p}.
\label{eq:M1-upper}
\end{equation}
If $E=T_{X_0}\cR_r$ for some $X_0\in\cR_r$, then
\begin{equation}
I_{T_r}(E)=\frac{m_rc_q}{p}.
\label{eq:M1-equality}
\end{equation}
Consequently, every tangent plane of $\cR_r$ attains the global maximum in \eqref{eq:moderation-one}.
\end{proposition}

\begin{proof}
Let $Q=W^\perp$ and choose a Frobenius-orthonormal basis $A,B$ of $Q$. For $([v],K)$ in the range of the Gauss map, choose a unit representative $v$ and an orthonormal basis $k_1,k_2$ of $K$. By Lemmas~\ref{lem:gauss-haar} and~\ref{lem:hodge},
\begin{align}
I_{T_r}(W)&=m_r\int_{\mathbb{RP}^{p-1}}\int_{\Gr_2(\R^q)}\abs{\ip{A\wedge B}{(vk_1^T)\wedge(vk_2^T)}}\,d\nu_q(K)\,d\sigma_p([v]).
\label{eq:I-as-normal-average}
\end{align}
By the basic calculus of exterior algebra:\[
\langle x_1\wedge x_2,y_1\wedge y_2\rangle=\det\begin{pmatrix}
\langle x_1,y_1\rangle&
\langle x_1,y_2\rangle\\
\langle x_2,y_1\rangle&
\langle x_2,y_2\rangle
\end{pmatrix}\]
and 
\[
\ip{A}{vk^T}_F=\ip{A^Tv}{k}_{\R^q},
\]
we have:
\begin{align}
\ip{A\wedge B}{(vk_1^T)\wedge(vk_2^T)}=\det\begin{pmatrix}
\ip{A^Tv}{k_1}&\ip{A^Tv}{k_2}\\
\ip{B^Tv}{k_1}&\ip{B^Tv}{k_2}
\end{pmatrix}=\ip{A^Tv\wedge B^Tv}{k_1\wedge k_2}.
\label{eq:wedge-identity}
\end{align}
The two-vector $A^Tv\wedge B^Tv$ is decomposable. By the $O(q)$-invariance of $c_q$,
\begin{equation}
\int_{\Gr_2(\R^q)}\abs{\ip{A^Tv\wedge B^Tv}{k_1\wedge k_2}}\,d\nu_q(K)=c_q\abs{A^Tv\wedge B^Tv}.
\label{eq:K-average}
\end{equation}
Consequently substituting \eqref{eq:K-average} and \eqref{eq:wedge-identity} into \eqref{eq:I-as-normal-average}:
\begin{equation}
I_{T_r}(W)=m_rc_q\int_{\mathbb{RP}^{p-1}}\abs{A^Tv\wedge B^Tv}\,d\sigma_p([v]).
\label{eq:I-v-average}
\end{equation}
The integrand is even in $v$, so we may use normalized spherical measure $\widetilde\sigma_p$ on $\mathbb S^{p-1}$. The Cauchy--Schwarz inequality, together with the identity $\int_{\mathbb S^{p-1}}vv^T\,d\widetilde\sigma_p(v)=\frac1pI_p$ give:
\begin{align*}
\int\abs{A^Tv\wedge B^Tv}\,d\widetilde\sigma_p(v)
&\leq \int\abs{A^Tv}\,\abs{B^Tv}d\widetilde\sigma_p(v)\\
&\leq\left(\int\abs{A^Tv}^2d\widetilde\sigma_p(v)\right)^{1/2}\left(\int\abs{B^Tv}^2d\widetilde\sigma_p(v)\right)^{1/2}\\
&=\frac{\norm{A}_F\norm{B}_F}{p}=\frac1p.
\end{align*}
This proves \eqref{eq:M1-upper}.

We next prove the second assertion. Let $E=T_{X_0}\cR_r$ and write its normal plane as $E^\perp=v_0\otimes K_0$. Choose an orthonormal basis $k_1^0,k_2^0$ of $K_0$ and set
\[
A=v_0(k_1^0)^T,\qquad B=v_0(k_2^0)^T.
\]
Then $A,B$ are Frobenius-orthonormal basis of $E^\perp$. For every unit $v$,
\[
A^Tv=\ip{v_0}{v}k_1^0,\qquad B^Tv=\ip{v_0}{v}k_2^0,
\]
and hence $\abs{A^Tv\wedge B^Tv}=\abs{\ip{v_0}{v}}^2. $ The normalized spherical average is $1/p$, substitution in \eqref{eq:I-v-average} proves \eqref{eq:M1-equality}.
\end{proof}

Proposition \ref{prop:M1} proves the first moderation condition. It remains to show the second.

\begin{proposition}
\label{prop:M2}
Let $E=T_{X_0}\cR_r$ be any tangent plane of $\cR_r$. Then
\begin{equation}
\bM((\pi_E)_\#T_r)=\int\abs{\ip{E}{\operatorname{Tan}(T_r,X)}}\,d\norm{T_r}(X)=I_{T_r}(E)>0.
\label{eq:M2}
\end{equation}
\end{proposition}

\begin{proof}
The group $O(p)\times O(q)$ acts transitively on the collection of normal planes $v\otimes K$, it therefore suffices to prove the assertion for $E_0=Q_0^\perp$, where $Q_0=e_p\otimes\operatorname{span}\{e_{r+1},e_{r+2}\}.$

Write
\[
X=\begin{pmatrix}B&C\\D&F\end{pmatrix},\qquad\text{where}\quad\begin{aligned}&B\in M_{r,r}(\R),\quad C\in M_{r,2}(\R),\\& D\in M_{1,r}(\R),\quad F\in M_{1,2}(\R). 
\end{aligned}
\]
Then
\[
E_0=\left\{\begin{pmatrix}B&C\\D&0\end{pmatrix}:B,C,D\text{ arbitrary}\right\},
\]
and $\pi_{E_0}$ forgets the block $F$.

Fix
\[
y=\begin{pmatrix}B&C\\D&0\end{pmatrix}\in E_0
\]
with $\det B\neq0$. Every element in the fiber over $y$ has the form
\[
X_F=\begin{pmatrix}B&C\\D&F\end{pmatrix}.
\]
Schur complement gives
\begin{equation}
\rank X_F=r+\rank(F-DB^{-1}C).\label{eq:schur-rank}
\end{equation}
Hence
\begin{equation}
X_F\in\cC_r\quad\Longleftrightarrow\quad F=DB^{-1}C.\label{eq:unique-completion}
\end{equation}
Thus the fiber of $\pi_{E_0}$ contains at most one point whenever $\det B\neq0$. 

The exceptional set
\[
Z=\left\{\begin{pmatrix}B&C\\D&0\end{pmatrix}\in E_0:\det B=0\right\}
\]
is the zero set of a nonzero polynomial on the $d$-dimensional space $E_0$, so $\cH^d(Z)=0$. 

Set $M=\cR_r\cap \mathbb B_1(0)$. By \eqref{eq:sing-dimension} and the fact that $\cC_r\cap \mathbb S^{pq-1}$ has dimension at most $d-1$, $M$ is a multiplicity-one rectifiable representative of $T_r$ up to a $\cH^d$-null set. Let $J(X)=J_d(\pi_{E_0}|_{T_XM})$ be the tangential Jacobian and define the effective fiber count by
\[
N(y)=\#\{X\in M:\pi_{E_0}(X)=y,\ J(X)>0\}.
\]
The area formula, applied to the sets where $J$ is positive and zero, gives
\begin{equation}
\int_MJ(X)\,d\cH^d(X)=\int_{E_0}N(y)\,d\cH^d(y).
\label{eq:area-formula}
\end{equation}
The rectifiable push-forward theorem identifies the multiplicity $\theta$ of $(\pi_{E_0})_\#T_r$ by
\begin{equation}
\theta(y)=\sum_{\substack{X\in M,\ \pi_{E_0}(X)=y\\J(X)>0}}1\pmod 2\qquad\text{for $\cH^d$-almost every }y\in E_0.
\label{eq:pushforward-multiplicity}
\end{equation}
By \eqref{eq:unique-completion}, $N(y)\in\{0,1\}$ outside $Z$, and hence $\theta(y)=N(y)$ for almost every $y\in E_0\setminus Z$. Since $\cH^d(Z)=0$, we have:
\begin{equation}
\bM((\pi_{E_0})_\#T_r)=\int_{E_0}\abs{\theta(y)}\,d\cH^d(y)=\int_{E_0}N(y)\,d\cH^d(y).
\label{eq:no-cancellation}
\end{equation}

For the orthogonal projection $\pi_{E_0}$, $J(X)=\abs{\ip{E_0}{T_XM}}.$ Combining \eqref{eq:area-formula}, \eqref{eq:no-cancellation} and Proposition~\ref{prop:M1} gives
\[
\bM((\pi_{E_0})_\#T_r)=I_{T_r}(E_0)=\frac{m_rc_q}{p}>0.
\]
This completes the proof. 
\end{proof}

The two preceding propositions now lead to the required moderation statement:

\begin{proposition}
\label{prop:determinantal-moderation}
The current $T_r$ is a moderation in the sense of Definition~\ref{def:moderation}, and
\[
\MOD(T_r)=\frac{m_rc_q}{p}>0.
\]
\end{proposition}

\begin{proof}
By \eqref{eq:sing-dimension}, the approximate tangent planes of $T_r$ are the tangent planes of $\cR_r$ for $\norm{T_r}$-almost every point. Proposition~\ref{prop:M1} verifies condition \eqref{eq:moderation-one} for every such plane, while Proposition~\ref{prop:M2} verifies condition \eqref{eq:moderation-two}. The common value is $m_rc_q/p>0$.
\end{proof}

\subsection{\texorpdfstring{Mod-$2$ minimization of $C(r+1,r+2,r)$}{Mod-2 minimization of C(r+1,r+2,r)}}

The previous section verified all hypotheses we need, we may now apply Liu's criterion:
\begin{theorem}\label{thm:main}
For every $r\ge1$, the current $T_r$ is mass-minimizing among all finite-mass integral currents $S\in\mathbf I_d(V;\Z_2)$ with $\partial S=\partial T_r$. Consequently, $\cC_r$ is an area-minimizing cone modulo $2$. \end{theorem}
\begin{proof}
By Proposition~\ref{prop:determinantal-moderation}, $T_r$ is a moderation. Theorem~\ref{thm:liu}, therefore implies that $T_r$ is mass-minimizing among finite-mass integral currents modulo $2$ whose difference from $T_r$ is the boundary of a finite-mass $(d+1)$-current.

Now let $S\in\mathbf I_d(V;\Z_2)$ satisfy $\partial S=\partial T_r$. Then $S-T_r$ is a finite-mass integral cycle modulo $2$ in the Euclidean space $V$. The Federer--Fleming isoperimetric theorem gives $Q\in\mathbf I_{d+1}(V;\Z_2)$ with $S-T_r=\partial Q.$ Applying Liu's criterion gives $\bM(T_r)\leq\bM(S)$ for every such $S\in\mathbf I_d(V;\Z_2)$. Moreover, after a rescaling, the same conclusion holds for every truncation
\[
T_{r,R}:=\llbracket\cC_r\rrbracket_2\mathbin{\llcorner}\mathbb B_R(0),\qquad R>0.
\]
Thus every truncated cone is minimizing for its spherical boundary. This implies the local statement for the full cone. 

Now let $C_r^\infty:=\llbracket\cC_r\rrbracket_2$ and let $S\in\mathbf I_{d,\mathrm{loc}}(V;\Z_2)$ satisfy
\[
\partial S=\partial C_r^\infty=0,\qquad K:=\spt(S-C_r^\infty)\text{ is compact}.
\]
Choose a slicing-good radius $R$ with $K\Subset \mathbb B_R(0)$. The two currents agree near $\partial\mathbb B_R(0)$, the minimizing property of $T_{r,R}$ then yields
\[
\bM(C_r^\infty\mathbin{\llcorner}\mathbb B_R(0))\leq\bM(S\mathbin{\llcorner}\mathbb B_R(0)).
\]
This proves the second assertion.
\end{proof}

For the case $r=1$, by the identification $C(2,3,1)\cong C(\mathbb S^1\otimes\mathbb S^2)$, Theorem \ref{thm:main} also gives an affirmative answer to \cite[Question 15]{Pedersen93}.

\section{Instability of \texorpdfstring{$C(p,p,p-1)$}{C(p,p,p-1)}}\label{sec:unstability-Cppp-1}

The regular strata of real determinantal varieties are minimal submanifolds of the corresponding Euclidean matrix spaces; see \cite{DET1,DET2}. Kerckhove--Lawlor proved that determinantal cones $C(p,p,p-1)$ are not area-minimizing \cite[Theorem 8.3(1)]{KerckhoveLawlor1999}. We now further prove that this family of cones are in fact unstable for every $p\ge 2$.

\begin{proposition}
\label{prop:S1tensorS1-unstable}
The minimal cone $\mathcal C=C(\mathbb S^{1}\otimes \mathbb S^{1})\subset\mathbb R^{2}\otimes\mathbb R^{2}\cong \mathbb R^{4}$ is unstable.
\end{proposition}
\begin{proof}
    Under the identification $\mathbb R^2\otimes\mathbb R^2\simeq M_{2\times2}(\mathbb R),$ the cone $C(\mathbb S^1\otimes \mathbb S^1)$ the determinant cone $\{X:\det X=0\}$. By an orthogonal change\[
y_1=\frac{x_1+x_4}{\sqrt2},\quad
y_2=\frac{x_2-x_3}{\sqrt2},\quad
y_3=\frac{x_1-x_4}{\sqrt2},\quad
y_4=\frac{x_2+x_3}{\sqrt2},
    \]
this becomes $y_1^2+y_2^2=y_3^2+y_4^2.$ Hence $C(\mathbb S^1\otimes \mathbb S^1)$ is orthogonally equivalent to the Clifford cone $C(\mathbb S^1\times \mathbb S^1)$, and is unstable by \cite{alm66}; see also \cite[\S 6.1]{JS68}.
\end{proof}

By \cite[Proposition~4.2]{KerckhoveLawlor1999} and \eqref{eq:identicication}, the cone $C(\mathbb S^1\otimes \mathbb S^1)$ is the local model along the rank-$p-2$ stratum of every $C(p,p,p-1)$. This observation leads to the instability of the whole family.

\begin{proposition}
  The determinantal cones $C(p,p,p-1),\,p\ge 2$ are all unstable.
\end{proposition}
\begin{proof}
The case $p=2$ is Proposition~\ref{prop:S1tensorS1-unstable}. We therefore assume $p\geq 3$, and write
\[
\mathcal D_p:=C(p,p,p-1)=\{X\in M_p(\mathbb R):\det X=0\}.
\]
By Proposition~\ref{prop:S1tensorS1-unstable}, the cone
$C(2,2,1)\subset M_2(\mathbb R)$ is unstable. Hence there exists
\[
\varphi\in C_c^\infty(C(2,2,1)_{\mathrm{reg}})
\]
such that
\[
Q_{C(2,2,1)}(\varphi):=\int_{C(2,2,1)_{\mathrm{reg}}}\left(|\nabla \varphi|^2-|A_{C(2,2,1)}|^2\varphi^2\right)\,d\mu_{C(2,2,1)}<0.
\]

We first observe that every Euclidean cylinder over $C(2,2,1)$ is unstable. Set $m:=p^2-4.$ Choose $\eta\in C_c^\infty(\mathbb R^m)$ such that $\eta\equiv1$ on $\mathbb B_1(0)$, and define
\[
\eta_R(y):=\eta(y/R),\qquad u_R(y,z):=\eta_R(y)\varphi(z).
\]
Since the second fundamental form of $\mathbb R^m\times C(2,2,1)$ is entirely contributed by the $C(2,2,1)$-factor, we have $|A_{\mathbb R^m\times C(2,2,1)}|^2(y,z)=|A_{C(2,2,1)}|^2(z).$ Consequently,
\begin{align*}
Q_{\mathbb R^m\times C(2,2,1)}(u_R)&=
\left(\int_{\mathbb R^m}\eta_R^2\,dy\right)Q_{C(2,2,1)}(\varphi)+
\left(\int_{\mathbb R^m}|\nabla\eta_R|^2\,dy\right)
\left(\int_{C(2,2,1)_{\mathrm{reg}}}\varphi^2\,d\mu_{C(2,2,1)}\right)\\
&=R^m\left(\int_{\mathbb R^m}\eta^2\,dy\right)Q_{C(2,2,1)}(\varphi)\\
&\quad+R^{m-2}
\left(\int_{\mathbb R^m}|\nabla\eta|^2\,dy\right)
\left(\int_{C(2,2,1)_{\mathrm{reg}}}\varphi^2\,d\mu_{C(2,2,1)}\right).
\end{align*}
Since $Q_{C(2,2,1)}(\varphi)<0$, we may choose $R$ sufficiently large so that $Q_{\mathbb R^m\times C(2,2,1)}(u_R)<0$. Fix such an $R$, and write $u=u_R$. It remains to transfer this negativity from the cylinder to the cone. We do so by blowing up at a point with rank $ p-2.$

Choose $X_0\in\mathcal D_p$ with $\operatorname{rank}X_0=p-2$. By \cite[Proposition~4.2]{KerckhoveLawlor1999}, we have the following local model for tangent cone:
\[
\operatorname{Tan}_{X_0}\mathcal D_p\cong\mathbb R^{p^2-4}\times C(2,2,1).
\]
Indeed, the local representation in \cite[Proof of Proposition~4.2]{KerckhoveLawlor1999} is
\[
\begin{pmatrix}
M+\rho A&\rho B\\
\rho C&\rho D+\rho^2 C(M+\rho A)^{-1}B
\end{pmatrix},
\]
after translation by $X_0$ and dilation by $\rho^{-1}$, the smooth error term is $\rho C(M+\rho A)^{-1}B.$ Consequently, letting $\mathcal D_{p,\rho}:=\rho^{-1}(\mathcal D_p-X_0),$ then
\[
\mathcal D_{p,\rho}^{\mathrm{reg}}\xrightarrow{C^\infty}\mathbb R^{p^2-4}\times C(2,2,1)_{\mathrm{reg}}\qquad\text{as }\rho\to 0
\]
on compact subsets of the regular part.

Choose $U\Subset\mathbb R^{p^2-4}\times C(2,2,1)_{\mathrm{reg}}$ with $\operatorname{spt}u\subset U$. By the smooth convergence of $\mathcal D_{p,\rho}^{\mathrm{reg}}$, for all sufficiently small $\rho$, there exist diffeomorphisms $\Phi_\rho:U\to \mathcal D_{p,\rho}^{\mathrm{reg}}$ onto their images such that $\Phi_\rho\to\mathrm{Id}$ in $C^\infty$. Define $u_\rho:=u\circ\Phi_\rho^{-1}.$ Then $u_\rho$ is compactly supported, and the smooth convergence therefore gives
\[
Q_{\mathcal D_{p,\rho}}(u_\rho)\longrightarrow
Q_{\mathbb R^{p^2-4}\times C(2,2,1)_{\mathrm{reg}}}(u)<0.
\]
Hence $Q_{\mathcal D_{p,\rho}}(u_\rho)<0$ for all sufficiently small $\rho$, so $\mathcal D_{p,\rho}$ is unstable. Since stability is invariant under Euclidean translations and dilations, the instability of $\mathcal D_{p,\rho}$ implies that of $\mathcal D_p$.
Therefore $C(p,p,p-1)$ is unstable for every $p\ge2$.
\end{proof}

\section{Full classification of Segre cones}\label{sec:Segre}

Very recently, the author Cui show that all Segre cones\[C\bigl(\mathbb S^{n_1-1}\otimes\cdots\otimes\mathbb S^{n_d-1}\bigr)\hookrightarrow\mathbb R^{\,n_1\cdots n_d}\] with $d\ge 3$ and dimension greater than $4$ are all regular area-minimizing cones \cite{CUI26}. The case $d=1$ is trivial; case $d=2$ has been fully classified by \cite{Pedersen93,KerckhoveLawlor1999}, Sections \ref{sec:remaining-real-determinantal} and \ref{sec:unstability-Cppp-1}, and the identification \eqref{eq:identicication}. Therefore the only exceptional case remained is the triple tensor cone $C(\mathbb S^1\otimes\mathbb S^1\otimes\mathbb S^1)$. The main difficulty is that the vanishing angle does not exist for the cone, and Lawlor's curvature criterion in \cite{Lawlorcriterion} does not apply.

In this section, we identify the triple tensor cone with a Harvey--Lawson calibrated cone constructed in \cite[\S 3, Theorem 3.1]{HL82}. This completes the full classification of all Segre cones as in the following Table \ref{tab:segre-cones}:

\begin{table}[H]
\centering
\small
\renewcommand{\arraystretch}{1.25}
\begin{tabularx}{\textwidth}{
  >{\centering\arraybackslash}p{0.12\textwidth}
  >{\centering\arraybackslash}p{0.25\textwidth}
  >{\centering\arraybackslash}p{0.16\textwidth}
  >{\raggedright\arraybackslash}X}
\toprule
$d$ & Segre cone & Dimension & Status and Source \\
\midrule
$d=1$
& $C(\mathbb S^{n_1-1})=\mathbb R^{n_1}$
& $n_1$
& Trivial \\
\addlinespace[0.6ex]
\multirow{3}{*}{$d=2$}
& $C(\mathbb S^1\otimes\mathbb S^1)$
& $3$
& Non-minimizing by \cite{KerckhoveLawlor1999}; in fact, unstable by Proposition~\ref{prop:S1tensorS1-unstable}. \\
& $C(\mathbb S^1\otimes\mathbb S^2)$
& $4$
& Mod 2-minimizing by Section~\ref{sec:remaining-real-determinantal}. \\
& $C(\mathbb S^{n_1-1}\otimes\mathbb S^{n_2-1})$, $n_1+n_2\ge6$
& $\ge5$
& Area-minimizing or mod 2-minimizing by \cite{Pedersen93,KerckhoveLawlor1999}. \\
\addlinespace[0.6ex]
\multirow{2}{*}{$d\ge3$}
& $C(\mathbb S^1\otimes\mathbb S^1\otimes\mathbb S^1)$
& $4$
& Area-minimizing by Proposition~\ref{prop:Segre-HL}. \\
& All other Segre cones
& $\ge5$
& Area-minimizing by \cite{CUI26}. \\
\bottomrule
\end{tabularx}
\caption{Classification of area-minimizing Segre cones.}
\label{tab:segre-cones}
\end{table}

We prove the assertion by explicitly constructing an orthogonal isomorphism between the triple tensor cone $C(\mathbb S^1\otimes \mathbb S^1\otimes \mathbb S^1)$ and the Harvey--Lawson cone as in the following Proposition \ref{prop:Segre-HL}:

\begin{proposition}
\label{prop:Segre-HL}
Let $V:=\mathbb R^2\otimes\mathbb R^2\otimes\mathbb R^2\cong \mathbb R^8$ and let $u(t):=(\cos t,\sin t)\in \mathbb S^1\subset\mathbb R^2.$ Consider the cone
\[
C:=C(\mathbb S^1\otimes \mathbb S^1\otimes \mathbb S^1)=\left\{r\,u(a)\otimes u(b)\otimes u(c):r\geq 0,\ a,b,c\in\mathbb R\right\}.
\]
Define the Harvey--Lawson cone in $\mathbb C^4\cong\mathbb R^8$ by
\[
C_{\mathrm{HL}}^4:=\left\{(z_1,z_2,z_3,z_4)\in\mathbb C^4:|z_1|=|z_2|=|z_3|=|z_4|,\quad iz_1z_2z_3z_4\in[0,\infty)\right\}.
\]
Then there exists a real orthogonal isomorphism $Q:V\to\mathbb C^4$ such that $Q(C)=C_{\mathrm{HL}}^4.$ In particular, $C$ is area-minimizing since $C_{\mathrm{HL}}^4$ is special Lagrangian and hence is calibrated and minimizing \cite[\S 3, Theorem 3.1]{HL82}.
\end{proposition}

\begin{proof}
Let $e_0=(1,0)$ and $e_1=(0,1)$ be the standard basis of $\mathbb R^2$. For $\varepsilon=(\varepsilon_1,\varepsilon_2,\varepsilon_3)\in\{0,1\}^3,$ 
write
\[
e_\varepsilon:=e_{\varepsilon_1}\otimes e_{\varepsilon_2}\otimes e_{\varepsilon_3},
\qquad|\varepsilon|:=\varepsilon_1+\varepsilon_2+\varepsilon_3.
\]
Thus every $x\in V$ can be written uniquely as $x=\sum_{\varepsilon\in\{0,1\}^3}x_\varepsilon e_\varepsilon.$

\medskip\noindent

\textbf{Step 1. Construction of the real orthogonal map.} Let
\begin{align*}
  \mathcal S_-&:=\left\{\sigma=(\sigma_1,\sigma_2,\sigma_3)\in\{\pm1\}^3:\sigma_1\sigma_2\sigma_3=-1\right\}\\
  &=\{(-1,-1,-1),(1,1,-1),(1,-1,1),(-1,1,1)\}.
\end{align*}
For $\sigma\in\mathcal S_-$, set $\sigma^\varepsilon:=
\sigma_1^{\varepsilon_1}
\sigma_2^{\varepsilon_2}
\sigma_3^{\varepsilon_3},$ and define $W_\sigma(x):=\frac12\sum_{\varepsilon\in\{0,1\}^3}i^{|\varepsilon|}\sigma^\varepsilon x_\varepsilon.$ We now define the real-linear map $Q:V\longrightarrow\mathbb C^4$ by
\[
Q(x):=e^{-i\pi/8}\bigl(
W_{(-1,-1,-1)}(x),
W_{(1,1,-1)}(x),
W_{(1,-1,1)}(x),
W_{(-1,1,1)}(x)\bigr).
\]

We verify that $Q$ is a real isometry. Recall that the real Euclidean inner product on $\mathbb C^4$ is $\langle z,w\rangle_{\mathbb R^8}=\operatorname{Re}\sum_{j=1}^4z_j\overline{w_j}.$ For $\delta\in\{0,1\}^3$, character orthogonality on $\mathcal S_-$ gives
\[
\sum_{\sigma\in\mathcal S_-}\sigma^\delta=
\begin{cases}
4, & \delta=(0,0,0),\\
-4, & \delta=(1,1,1),\\
0, & \text{otherwise}.
\end{cases}
\]
Consequently, for $x,y\in V$,
\begin{align*}
\langle Qx,Qy\rangle_{\mathbb R^8}=\operatorname{Re}\sum_{\sigma\in\mathcal S_-}W_\sigma(x)\overline{W_\sigma(y)}=\frac14\operatorname{Re}\sum_{\varepsilon,\eta\in\{0,1\}^3}i^{|\varepsilon|}(-i)^{|\eta|}x_\varepsilon y_\eta\sum_{\sigma\in\mathcal S_-}\sigma^{\varepsilon+\eta}.
\end{align*}
Here the exponents of the $\sigma_j$'s are understood modulo $2$. The diagonal terms $\eta=\varepsilon$ contribute $\sum_{\varepsilon\in\{0,1\}^3}x_\varepsilon y_\varepsilon.$ The only remaining potentially nonzero terms occur when $\eta=(1,1,1)-\varepsilon.$ Denote $k=|\varepsilon|$, then\[
i^{|\varepsilon|}(-i)^{|\eta|}=i^k(-i)^{3-k}\in i\mathbb R,
\]so these terms have vanishing real part. It follows that $\langle Qx,Qy\rangle_{\mathbb R^8}=\sum_{\varepsilon\in\{0,1\}^3}x_\varepsilon y_\varepsilon=\langle x,y\rangle_V.$ Thus $Q$ is a real orthogonal isomorphism.

\medskip\noindent

\textbf{Step 2. The inclusion $Q(C)\subset C_{\mathrm{HL}}^4$.} Let $x=r\,u(a)\otimes u(b)\otimes u(c)\in C.$ Writing $(t_1,t_2,t_3):=(a,b,c),$ its coordinates are $x_\varepsilon=r\prod_{j=1}^3(\cos t_j)^{1-\varepsilon_j}(\sin t_j)^{\varepsilon_j}.$ Therefore, for every $\sigma\in\mathcal S_-$,
\begin{align*}
W_\sigma(x)&=\frac r2\sum_{\varepsilon\in\{0,1\}^3}\prod_{j=1}^3
\left(
(i\sigma_j)^{\varepsilon_j}
(\cos t_j)^{1-\varepsilon_j}
(\sin t_j)^{\varepsilon_j}
\right)\\
&=\frac r2\prod_{j=1}^3\bigl(\cos t_j+i\sigma_j\sin t_j\bigr)=\frac r2\exp\left(i\sum_{j=1}^3\sigma_jt_j\right).
\end{align*}
It follows that
\begin{align}
Q(x)=\frac r2e^{-i\pi/8}\bigl(e^{-i(a+b+c)},e^{i(a+b-c)},e^{i(a-b+c)},e^{i(-a+b+c)}\bigr).
\label{eq:Q-rank-one}
\end{align}
In particular, if we denote $Q(x)=(z_1,z_2,z_3,z_4),$ then $|z_i|=\frac r2,$ $i=1,\ldots,4.$ Moreover, the four $a,b,c$-dependent phases in \eqref{eq:Q-rank-one} sum to zero. Hence
\begin{align*}
z_1z_2z_3z_4=\left(\frac r2\right)^4 e^{-4i\pi/8}=-i\left(\frac r2\right)^4.
\end{align*}
Thus $iz_1z_2z_3z_4=\left(\frac r2\right)^4\in[0,\infty)$, and therefore $Q(C)\subset C_{\mathrm{HL}}^4.$

\medskip\noindent

\textbf{Step 3. The reverse inclusion $C_{\mathrm{HL}}^4\subset Q(C)$.} Let $z=(z_1,z_2,z_3,z_4)\in C_{\mathrm{HL}}^4.$ If $z=0$, then $z=Q(0)$, so suppose that $z\neq0$. There exist $s>0$ and $\theta_1,\theta_2,\theta_3,\theta_4\in\mathbb R$ such that $z_j=se^{i\theta_j}$, $j=1,\cdots,4.$ Since $iz_1z_2z_3z_4>0$, we have
\[
\theta_1+\theta_2+\theta_3+\theta_4\equiv -\frac{\pi}{2}\pmod{2\pi}.
\]
Set $\varphi_j\equiv\theta_j+\frac{\pi}{8}\pmod{2\pi}.$ Then
\[
\varphi_1+\varphi_2+\varphi_3+\varphi_4\equiv0\pmod{2\pi}.
\]
Choose $\varphi_j\in\mathbb R$, such that $\varphi_1+\varphi_2+\varphi_3+\varphi_4=0.$ Define
\[
a:=\frac{\varphi_2+\varphi_3}{2},\qquad
b:=\frac{\varphi_2+\varphi_4}{2},\qquad
c:=\frac{\varphi_3+\varphi_4}{2}.
\]
Using $\varphi_1+\varphi_2+\varphi_3+\varphi_4=0,$ a direct calculation gives
\[
a+b-c=\varphi_2,\qquad
a-b+c=\varphi_3,\qquad
-a+b+c=\varphi_4,\qquad
-a-b-c=\varphi_1.
\]
Therefore
\[
\begin{pmatrix}
-(a+b+c)\\
a+b-c\\
a-b+c\\
-a+b+c
\end{pmatrix}=
\begin{pmatrix}
\varphi_1\\
\varphi_2\\
\varphi_3\\
\varphi_4
\end{pmatrix}.
\]
Taking $r:=2s$, formula \eqref{eq:Q-rank-one} yields
\begin{align*}
Q\bigl(r\,u(a)\otimes u(b)\otimes u(c)\bigr)&=s e^{-i\pi/8}
\bigl(
e^{i\varphi_1},
e^{i\varphi_2},
e^{i\varphi_3},
e^{i\varphi_4}
\bigr)
\\&=s\bigl(
e^{i\theta_1},
e^{i\theta_2},
e^{i\theta_3},
e^{i\theta_4}
\bigr)=z.
\end{align*}
Thus $C_{\mathrm{HL}}^4\subset Q(C)$ and hence $Q(C)=C_{\mathrm{HL}}^4$ as desired.
\end{proof}

Together with the results of \cite{KerckhoveLawlor1999,CUI26,Pedersen93} and previous sections, this identification completes the classification summarized in Table \ref{tab:segre-cones}.

\medskip\noindent
\paragraph{Declaration of Generative AI Use.} During the preparation of this work, the authors used ChatGPT and DeepSeek to improve the grammar, clarity, quality of the manuscript's writing, and to assist in creating tables. After using these tools, the authors reviewed and edited the content as needed and take full responsibility for the content of the article.

\paragraph{Conflict of Interest}On behalf of all authors, the corresponding author states that there is no conflict of interest.

\paragraph{Data Availability}Data sharing is not applicable to this article, as no datasets were generated or analyzed during the current study.

\paragraph{Funding} The author Cui is supported by the National Natural Science Foundation of China (Grant no. 12301068), the project of Stable Support for Youth Team in Basic Research Field, CAS (YSBR-001), and Xiaomi Young Scholar Fund. The authors Jiao and Zhu are supported by the National Natural Science Foundation of China (Grant no. 12371055).

\begin{flushleft}
\medskip\noindent
\begin{tabbing}

				Hongbin Cui\\
				School of Mathematical Sciences, University of Science and Technology of China\\
				Wu Wen-Tsun Key Laboratory of Mathematics, USTC, Chinese Academy of Sciences\\
				96 Jinzhai Road, Hefei, 230026, Anhui Province, China\\			
				\texttt{cuihongbin@ustc.edu.cn}
								
\end{tabbing}

\begin{tabbing}

				Xiaoxiang Jiao\\
				School of Mathematical Sciences, University of Chinese Academy of Sciences\\
				19A Yuquan Road, Beijing, 100049, China\\
\texttt{xxjiao@ucas.ac.cn}
				
\end{tabbing}

\begin{tabbing}
Hangyue Zhu\\School of Mathematical Sciences, University of Chinese Academy of Sciences\\
				19A Yuquan Road, Beijing, 100049, China\\
\texttt{zhuhangyue24@mails.ucas.ac.cn}
\end{tabbing}

\end{flushleft}

\bibliographystyle{alpha}
\bibliography{ref}

\end{document}